\documentclass[11pt,a4paper]{amsart}
\usepackage{amsthm}
\usepackage[utf8]{inputenc}
\usepackage[english]{babel}
\usepackage[T1]{fontenc}
\usepackage{url}
\usepackage{mathrsfs}
\usepackage{tikz}
\usepackage{upref}
\usepackage{cancel} 

\title{Neighboring mapping points theorem}

\author{Andrei V. Malyutin}
\thanks{The first author was supported by RFBR according to the research project n.~20-01-00070. The second author is supported by the Ministry of Science and Higher Education of the Russian Federation, agreement 075-15-2019-1620 date 08/11/2019.}
\address{A.\,V.~Malyutin: St.\,Petersburg Department of Steklov Institute of Mathematics, 
Fontanka 27, St. Petersburg 191011, Russia;
St.~Petersburg State University}
\email{malyutin@pdmi.ras.ru}

\author{Oleg R. Musin}
\address{O.\,R.~Musin: University of Texas Rio Grande Valley, School of Mathematical and
 Statistical Sciences, One West University Boulevard, Brownsville, TX, 78520; 
 Leonhard Euler International Mathematical Institute in St.\,Petersburg;
 Moscow Institute of Physics and Technology;
 The Institute for Information Transmission Problems of RAS}
\email{oleg.musin@utrgv.edu}

\newcommand \diamE   {\operatorname{diam_{Eu}}}
\newcommand \diamA   {\operatorname{diam_{A}}}
\newcommand \crA     {\operatorname{circ_{A}}}
\newcommand \sk      {\operatorname{sk}}
\newcommand \rk      {\operatorname{rk}}

\newcommand \R       {\mathbb R}
\newcommand \SSS     {\mathbb S}

\newcommand \dd      {\partial}

\newcommand \be      {\begin{equation}}
\newcommand \ee      {\end{equation}}

\newtheorem{thm}{Theorem}
\newtheorem{prop}{Proposition}
\newtheorem{lem}{Lemma}

\newtheorem{cor}{Corollary}

\theoremstyle{definition}
\newtheorem{defn}{Definition}

\theoremstyle{remark}
\newtheorem{rmrk}{Remark}

\newtheorem{example}{Example}

\begin{document}

\begin{abstract}  
We introduce and study a new family of theorems extending the class of Borsuk--Ulam and topological Radon type theorems. 
The defining idea for this new family is to replace requirements of the form 
`the image of a subset that is large in some sense is a singleton' 
with requirements of the milder form
`the image of a subset that is large in some sense is a subset that is small in some sense'.
This approach 
covers the case of mappings $\SSS^m\to\R^n$ with $m<n$ and extends to wider classes of spaces. 

An example of a statement from this new family is the following theorem.
Let $f$ be a continuous map of the boundary~$\partial\Delta^{n}$ of the $n$\nobreakdash-dimensional simplex~$\Delta^n$ to a contractible metric space~$M$.
Then $\partial\Delta^{n}$ contains a subset~$E$ such that $E$ (is `large' in the sense that~it) intersects all facets of $\Delta^n$ and the image $f(E)$ (is `small' in the sense that~it) is either a singleton or a subset of the boundary $\partial B$ of a metric ball~$B\subset M$ whose interior does not meet~$f(\partial\Delta^{n})$. 

We generalize this theorem to noncontractible normal spaces via covers and deduce a series of its corollaries. Several of these corollaries are similar to the topological Radon theorem.
\end{abstract} 

\maketitle

\section{Introduction}
In this paper we introduce and study a new family of theorems extending the class of Borsuk--Ulam and topological Radon type theorems (though none of our theorems is a generalization for the Borsuk--Ulam or topological Radon theorem itself). 
By the Borsuk--Ulam and topological Radon type theorems we mean those stating that a continuous map takes a `wide' set of some specific kind to a point. 
Let us list several of the most 
influential examples.

\begin{itemize}
\item {\it The Borsuk--Ulam theorem} itself says that every continuous map of a Euclidean $n$-sphere~$\SSS^n$ into Euclidean $n$-space~$\R^n$ identifies two antipodes.

\item {\it The Hopf theorem} states that if~$X$ is a compact Riemannian $n$-manifold and $f\colon X\to\R^n$ is a continuous map, then for any $\delta > 0$, there exists a geodesic $\gamma\colon[0,\delta]\to X$ of length~$\delta$ such that $f(\gamma(0))=f(\gamma(\delta))$.

\item {\it The topological Radon theorem} says that if $P$ is a convex $n$-polytope, then any continuous map $\partial P\to\R^{n-1}$ identifies two points from disjoint faces.

\item {\it The topological Tverberg theorem} says that if $d\ge 1$ is an integer, $r$ is a prime power, and $P$ is a convex $(r-1)(d+1)$-polytope, then any continuous map $\partial P\to \R^d$
identifies $r$ points from $r$ pairwise disjoint faces.
\end{itemize}

See, e.\,g., \cite{Ste85, Ste93, Mat03, Kar08, AKV12, Fri15, BBZ16, Sko18, BS18} and references therein for more examples including various extensions and generalizations for $\mathbb{Z}_p$-spaces, maps between manifolds, matroids, colored versions, etc. 
Another related family is Knaster's conjecture type theorems (see \cite{Mat14} and references therein).

All of these examples involve rigid dimensional restrictions.
It is a natural question whether 
the maps not satisfying these restrictions 
have any properties of the Borsuk--Ulam kind.
In particular, we are interested in whether the Borsuk--Ulam theorem has any reasonable extensions to the case of mappings $\SSS^m\to\R^n$ with $m<n$ (a related idea appears in~\cite{ABF20}).

Extensions of this kind are found in a new class we study.
This class emerges 
by replacing conditions of the form 
`\emph{the image of a subset that is large in some sense is a singleton}' 
with conditions of the milder form
`\emph{the image of a subset that is wide in some sense is a subset that is restricted in some sense}'.
This approach 
covers the case of mappings $\SSS^m\to\R^n$ with $m<n$ and extends to wider classes of spaces. 

Here is an example for the simplest nondegenerate case $\SSS^1\to\R^2$.

\begin{prop}[\cite{Mal16}]\label{prop:plane-curves}
Let $a$, $b$, and $c$ be three closed arcs covering the circle $\SSS^1$ such that no two of them cover $\SSS^1$,
and let $f\colon\SSS^1\to\R^2$ be a continuous map.
Then either $f(a)\cap f(b)\cap f(c)\neq\emptyset$ or each of $f(a)$, $f(b)$, and $f(c)$ touches a  closed Euclidean disk~$D^2\subset\R^2$ whose interior does not meet~$f(\SSS^1)$.
\end{prop}

\begin{figure}[h]
\scalebox{0.6}{
\begin{tikzpicture}
\draw[thick,rounded corners=20pt,densely dotted]
(0,-3) -- (0,0) -- (1,3)-- (5,4) -- (6,2) -- (1.45,1.35); 
\draw[thick,rounded corners=20pt]
(1.45,1.35) -- (-1,1) -- (-1, -1) -- (7,-2)--(8,-1);
\draw[thick,rounded corners=20pt,densely dashed]
(8,-1) -- (9,0) -- (8,5) -- (6,5) -- (4, -3) -- (0, -4) -- (0,-3);
\draw[gray] (3.16,0.03) circle (1.53) node[black]{\large $D^2$};
\draw (3.16,4.03)node{\large $f(a)$};
\draw (-1.66,0.03)node{\large $f(b)$};
\draw (9.06,3.03)node{\large $f(c)$};
\fill[gray] (4.66,-0.3) circle (0.07); 
\fill[gray] (3,1.55) circle (0.07); 
\fill[gray] (3,-1.5) circle (0.07); 
\draw[thick,densely dotted] (11,4)--(12,4) node[right]{\large $a$};
\draw[thick]                (11,3.5)--(12,3.5) node[right]{\large $b$};
\draw[thick,densely dashed] (11,3)--(12,3) node[right]{\large $c$};
\end{tikzpicture}}
\caption{For Proposition~\ref{prop:plane-curves}. The circle~$\partial D^2$ touches $f(a)$, $f(b)$, and $f(c)$}
\label{fig:touches}
\end{figure}

Proposition~\ref{prop:plane-curves} works for plane curves and knot diagrams and has a corollary with applications in knot theory (see~\cite{Mal16}). We formulate this corollary here.
Let $\gamma\colon \SSS^1\to\R^2$ be a regular smooth plane curve in general position (that is, its only singularities are transversal double points).
By an \emph{edge} of $\gamma$ we mean the closure of a component of the set $\gamma(\SSS^1)\setminus V$, where $V$ is the set of double points of~$\gamma$. We say that two edges~$I$ and~$J$ of~$\gamma$ are \emph{neighboring edges} or \emph{neighbors} if there exists a component~$Q$ of $\R^2 \setminus \gamma(\SSS^1)$ such that the boundary $\dd Q$ contains both $I$ and~$J$. We say that two edges~$I$ and~$J$ of~$\gamma$ are \emph{consecutive} if the union $I \cup J$ coincides with the image $\gamma(\alpha)$ of a (connected) arc $\alpha$ in~$\SSS^1$. We denote by $\rho$ the maximal metric on the set~$E(\gamma)$ of edges of $\gamma$ in the class of metrics satisfying the condition
`$\rho(I, J) = 1$ whenever $I$ and $J$ are consecutive edges of~$\gamma$'.

\begin{prop}[\cite{Mal16}]
\label{cor:curve}
If the curve~$\gamma$ has $k$ double points, then~$\gamma$ has a pair of
neighboring edges $I$ and $J$ with $\rho(I, J) \ge 2k/3$. 
\end{prop}

Proposition~\ref{prop:plane-curves} readily implies Proposition~\ref{cor:curve} if we choose the arcs $a$, $b$, and $c$ appropriately.
Proposition~\ref{cor:curve} appears in~\cite{Mal16} as an auxiliary lemma (Lemma~5.1) needed to obtain a series of statements related to knot theory.
In~\cite{Mal16} this lemma is deduced from the topological Helly theorem (see \cite{Bog02, Mon14}).  
The statement of Proposition~\ref{cor:curve} was one of the starting points for our study. 

In this paper we generalize Proposition~\ref{prop:plane-curves} to noncontractible normal spaces via covers.
The generalizations and their corollaries will be formulated in the next sections after definitions.
Our method is based on obstruction theory and uses a variation of the concept of non--null--homotopic covers introduced in~\cite{Mus16,Mus17}.

\section{Definitions and results}

Throughout this paper we mainly consider normal topological spaces\footnote{A topological space~$X$ is \emph{normal} if any two disjoint closed sets of~$X$ are contained in disjoint open sets of~$X$; see \cite[p.~446]{Sch97} for equivalent definitions via the Urysohn and shrinking lemmas.}, 
all simplicial complexes and covers will be finite, 
all manifolds will be both compact and PL, 
$\SSS^n$~will denote the $n$-dimensional sphere,
$\Delta^n$~will denote the $n$-dimensional simplex,
$\sk_k(\Delta^n)$ will denote the $k$-skeleton of~$\Delta^n$.
We shall denote the set of homotopy classes of continuous maps $V\to W$ by $[V,W]$. 
The nerve of a (finite) collection $\mathcal{S}$ of sets will be denoted by $\mathcal{N}(\mathcal{S})$.
When this does not cause confusion we use the same notation for an abstract simplicial complex and its underlying space (carrier).

The further exposition in this section is structured as follows: first we give a chain of successively stronger generalizations of Proposition~\ref{prop:plane-curves} (Theorem~\ref{thm:simplex} is the weakest, Theorem~\ref{thm:G-EPtriple} is the strongest one); then we present a family of corollaries (all but one of which follow from Theorem~\ref{thm:simplex}). 

\subsection{Spherical $f$--neighbors}
All of the following generalizations and corollaries replace the condition `the image is a singleton' appearing in the Borsuk--Ulam type theorems with the following milder condition of `spherical neighboring'. 

\begin{defn}[Spherical $f$-neighbors]
Let $X$ be a set, let $Y$ be a metric space, and
let ${f\colon X\to Y}$ be a map.
We say that a subset $N\subset X$ is a \emph{set of spherical $f$--neighbors} 
if $N$ contains at least two points and the image $f(N)$ is 
either a point or a subset of the boundary $\partial B$ of a metric ball\footnote{By a \emph{metric ball} in a metric space $(Y,d)$ with metric~$d$ we mean a subset of the form $\{y\in Y\mid d(y,x)\le R\}$,  $x\in Y$, $R\ge 0$.}~$B\subset Y$ whose interior does not meet~$f(X)$. 
If a two-point set $\{p,q\}$ is a set of spherical $f$--neighbors, we say that $p$ and $q$ are \emph{spherical $f$--neighbors}. 
(See Fig.~\ref{fig:neighbors}.)
\end{defn}

\begin{figure}[h]
\scalebox{0.6}{
\begin{tikzpicture}
\draw[thick,rounded corners=20pt]
(0,-1) -- (0,0) -- (1,3) -- (5,4) -- (6,2) -- (-1,1) -- (-1, -1) -- (7,-2) -- (9,0) -- (8,5) -- (6,5) -- (4, -3) -- (0, -4) -- (0,-1);
\draw[gray] (3.16,0.03) circle (1.53) node[black]{\large $B$};
\fill[gray] (4.66,-0.3) circle (0.1) node[right,black]{\large $f(p_2)$};
\fill[gray] (3,1.55) circle (0.1) node[above,black]{\large $f(p_1)$};
\fill[gray] (3,-1.5) circle (0.1) node[below,black]{\large $f(p_3)$};
\draw (-2.36,0.03)node{\large $f(X=\SSS^1)$};
\end{tikzpicture}}
\caption{The image of a set $\{p_1,p_2,p_3\}$ of spherical $f$--neighbors}
\label{fig:neighbors}
\end{figure}

The first extension generalizes Proposition~\ref{prop:plane-curves} to the case of spheres of arbitrary dimension and replaces Euclidean spaces with arbitrary contractible metric spaces. 
(Recall that \emph{facets} of a polytope of dimension $n$ are its faces of dimension~$n-1$.)

\begin{thm}\label{thm:simplex}
Let $f$ be a continuous map of the boundary $\dd\Delta^{n}$ of the $n$\nobreakdash-dimensional simplex $\Delta^n$ to a contractible metric space~$M$.
Then a set of spherical $f$--neighbors intersects all facets of~$\Delta^{n}$.
\end{thm}

\begin{proof}
Let $d$ denote the metric on~$M$. If $z$ is a point and $N$ is a subset in~$M$, we write
$$
d(z,N):=\inf_{p\in N}d(z,p).
$$
Let $\Delta_1, \dots, \Delta_{n+1}$ be the facets of~$\Delta^{n}$.
For each $i\in\{1,\dots,n+1\}$, we set
$$
E_i:=\{z\in M \mid d(z,f(\Delta_i))=d(z,f(\dd\Delta^{n}))\}.
$$ 
Observe that $E_i$ contains~$f(\Delta_i)$ and is closed,
so that $\{E_1,\dots,E_{n+1}\}$ is a closed cover of~$M$.
Since $M$ is contractible, $f$ extends to a continuous map~$F\colon \Delta^{n}\to M$.
Then $\{F^{-1}(E_1),\dots,F^{-1}(E_{n+1})\}$ is a closed cover of~$\Delta^{n}$ extending the closed cover $\{\Delta_1, \dots, \Delta_{n+1}\}$ of~$\dd\Delta^{n}$.
By the Knaster--Kuratowski--Mazurkiewicz (KKM) lemma the elements of $\{F^{-1}(E_1),\dots,F^{-1}(E_{n+1})\}$ have a common point~$p$. 
Then $x=F(p)$ lies in $E_1\cap\dots\cap E_{n+1}$.
Then either $x\in f(\dd\Delta^{n})$ so that $x$ belongs to all of~$f(\Delta_i)$ by the definition of~$E_i$, or $r=d(x,f(\dd\Delta^{n}))>0$ and the ball $B_r(x)$ of radius~$r$ centered at~$x$ touches all of~$f(\Delta_i)$ while its interior does not meet~$f(\partial\Delta^{n})$. 
This implies the statement.
\end{proof}

We generalize Theorem~\ref{thm:simplex} by replacing the set of facets with a more general class of covers as in the KKM lemma.

\subsection{KKM covers and spherical $f$--neighbors}

\begin{defn}[KKM covers] 
\label{def:KKM-covers}
Let $\Delta^{n+1}$ be an $(n+1)$-dimensional simplex with vertices labelled $v_1,\ldots,v_{n+2}$.
A~closed cover $\{C_1,\ldots,C_{n+2}\}$ of the $n$\nobreakdash-sphere~$\SSS^n$ is called a \emph{KKM cover} if there exists a homeomorphism $h\colon\SSS^n\to \partial \Delta^{n+1}$ such that for each $J\subset \{1,\ldots,n+2\}$ 
the convex hull of the vertices $v_j$ with $j\in J$ is covered by the union $\bigcup _{j\in J}h(C_j)$.
\end{defn}

The argument in the proof of Theorem~\ref{thm:simplex} also proves the following theorem. 

\begin{thm}\label{thm:KKM}
Let $\mathcal C$ 
be a KKM cover of the $n$\nobreakdash-sphere~$\SSS^n$,
and let $f\colon \SSS^n\to M$ be a continuous map to a contractible metric space~$M$.
Then a set of spherical $f$--neighbors intersects all elements of~$\mathcal C$.
\end{thm}

The key role in Theorem~\ref{thm:KKM} is played by the properties of the cover, and not by the fact that the underlying space is a sphere. 
To move on to the next generalization, 
we define non--null--homotopic covers (we generalize the concept of non--null--homotopic covers given in \cite{Mus16,Mus17}). 

\subsection{Non--null--homotopic covers and spherical $f$--neighbors}

Let $X$ be a normal topological space, and let $\mathcal U=\{U_1,\ldots,U_n\}$ be an open cover of~$X$.  
Let $\mathcal{N}(\mathcal{U})$ be the nerve of~$\mathcal U$.
Let $\Phi=\{\varphi_1,\ldots,\varphi_n\}$ be a partition of unity subordinate to $\mathcal U$. 
Let $v_1,\ldots,v_n$ be the vertices of the $(n-1)$-dimensional unit simplex~$\Delta^{n-1}$, where 
$$
\Delta^{n-1}:=\{x\in\R^n \mid x_i\ge0,~x_1 +...+x_n =1\}.
$$
For each $i$ we identify the vertex of $\mathcal{N}(\mathcal{U})$ corresponding to~$U_i$ with $v_i$ so that $\mathcal{N}(\mathcal{U})$ becomes a subcomplex of~$\Delta^{n-1}$.
We set
$$
h_{\mathcal U,\Phi}(x):=\sum\limits_{i=1}^n{\varphi_i(x)v_i}.
$$
Then $h_{\mathcal U,\Phi}$ is a continuous map from $X$ to $\mathcal{N}(\mathcal{U})\subset \Delta^{n-1}$. 
Since the linear homotopy $\Theta(t)=(1-t)\Phi + t\Psi$ of two partitions of unity $\Phi$ and $\Psi$ subordinate to $\mathcal U$ induces a homotopy between the corresponding maps,  
it follows that the homotopy class $[h_{\mathcal U,\Phi}]$ in $[X,\mathcal{N}(\mathcal{U})]$, 
where by $[V,W]$ we denote the set of homotopy classes of continuous maps $V\to W$,
does not depend on~$\Phi$ (see~\cite[Lemma 1.6]{Mus16}).  
We denote this class in $[X,\mathcal{N}(\mathcal{U})]$~by $[\mathcal U]$.  

The homotopy classes of covers are also well defined for closed sets. 
Indeed, in a normal space any finite closed cover has an open extension with the same nerve (see, e.\,g.,~\cite[Theorem~1.3]{Mor50} and \cite[pp.~31--33]{Iva20}).
Furthermore, if $\mathcal C=\{C_1,\ldots,C_n\}$ is a closed cover of a normal space~$X$ and $\mathcal S=\{S_1,\ldots,S_n\}$ and $\mathcal U=\{U_1,\ldots,U_n\}$ are two open covers 
such that $S_i\cap U_i$~contains~$C_i$ for all~$i$ and with the same nerve~$\mathcal{N}(\mathcal S)=\mathcal{N}(\mathcal U)=\mathcal{N}(\mathcal C)$, then each partition of unity subordinate to the open cover 
$$\mathcal T:=\{S_1\cap U_1,\ldots,S_n\cap U_n\}$$ is also subordinate to both $\mathcal S$ and $\mathcal U$. 
This implies that $[\mathcal S]=[\mathcal T]=[\mathcal U]$ in $[X,\mathcal{N}(\mathcal C)]$ due to the independence of the choice of partition of unity mentioned above.  
Then we set 
$$[\mathcal C]:=[\mathcal S]=[\mathcal T]=[\mathcal U].$$

\begin{defn}[Non--null--homotopic covers] 
\label{def2}
We say that an open or closed cover~$\mathcal C$ of a normal topological space~$X$ is {\em non--null--homotopic} if the corresponding homotopy class~$[\mathcal C]$ in $[X,\mathcal{N}(\mathcal C)]$ contains no constant map. 
\end{defn}

\begin{rmrk}
\label{rmrk:1}
Any non--null--homotopic map $X\to K$ to a finite simplicial complex yields non--null--homotopic covers on~$X$; to obtain an example, take the inverse images of all elements in one of the following collections:
\begin{itemize}
\item open stars of vertices of~$K$, 
\item stars of vertices of~$K$ in its first barycentric subdivision, 
\item maximal simplexes of~$K$.
\end{itemize}
\end{rmrk}

\begin{defn}[Homotopy ranks of maps] 
\label{def:homotopy-ranks}
Let $X$ be a topological space, let $K$ be a finite simplicial complex, and let $h\colon X\to K$ be a continuous map. Let $\Delta_K$ be the simplex spanned by the vertices of~$K$ so that $K$ is a subcomplex of~$\Delta_K$. We define the \emph{homotopy rank} $\rk(h)$ of~$h$ to be the least nonnegative integer~$k$ such that $h$ is null--homotopic in $K\cup\sk_k(\Delta_K)$, where $\sk_k$ stands for the $k$-skeleton.\footnote{We say that $K\cup\sk_k(\Delta_K)$ is the \emph{$k$-exoskeleton} of~$K$.}
(Since $\Delta_K$ is contractible, the homotopy rank is well defined and does not exceed the dimension of~$\Delta_K$.)
\end{defn}

\begin{rmrk}
In terms of Definition~\ref{def:homotopy-ranks}, $h$ is null--homotopic if and only if $\rk(h)=0$.
\end{rmrk}

\begin{defn}[Ranks of covers] 
We define the (\emph{homotopy}) \emph{rank} $\rk(\mathcal C)$ of a (closed or open) finite cover~$\mathcal C$ of a normal space~$X$ to be the homotopy rank of maps $X\to\mathcal{N}(\mathcal C)$ in the class~$[\mathcal C]$ determined by~$\mathcal C$.
\end{defn}

\begin{rmrk}
A cover is non--null--homotopic if and only if it is of nonzero rank.
\end{rmrk}

\begin{rmrk}
Since $\sk_m(\Delta^{m})=\Delta^{m}$ is contractible, it follows that the rank of an $n$-element cover does not exceed~$n-1$. 
\end{rmrk}

\begin{defn}[Principal covers] 
An $n$-element cover ($n\ge 2$) of rank $n-1$ is said to be \emph{principal}. 
\end{defn}

\begin{rmrk}
\label{rmrk:principal-criterion}
Since any proper nonempty subcomplex of~$\partial\Delta^m$ is contractible in~$\partial\Delta^m$, it follows that a cover is principal if and only if it is non--null--homotopic and its nerve is the boundary of a simplex. 
\end{rmrk}

\begin{rmrk}
\label{rmrk:principal-disjoint}
Remark~\ref{rmrk:principal-criterion} implies that no principal cover has a proper subcollection of elements with empty intersection;
in particular, no principal cover has disjoint elements.  
\end{rmrk}

\begin{rmrk}
Any non--null--homotopic map $X\to \SSS^k$ to the $k$-sphere yields a principal cover of~$X$ of rank $k+1$ (cf.~Remark~\ref{rmrk:1}; cf.~\cite[Theorem~1.5]{Mus16}). Thus, a space can have principal covers of distinct ranks.
\end{rmrk}

\begin{rmrk}[Conditions for cover non--null--homotopicity]
\label{rmrk:suf-cond}
If the composition of continuous maps is non--null--homotopic, then each of them is non--null--homotopic.
\begin{itemize}
\item On the one hand this implies that any refinement of a cover of rank~$k$ has rank at least~$k$. 
In particular, any refinement of a non--null--homotopic cover is non--null--homotopic. 

\item On the other hand this implies that if $f\colon X\to Y$ is a continuous map of normal spaces and $\mathcal{C}=\{C_1,\ldots,C_n\}$ is a closed cover of~$Y$ such that the dimension of the nerve~$\mathcal{N}(\mathcal{C})$ is less than the rank $\rk(\tilde{\mathcal{C}})$ of the induced cover $\tilde{\mathcal{C}}=\{f^{-1}(C_1),\ldots,f^{-1}(C_n)\}$, then $\rk(\mathcal{C})\ge \rk(\tilde{\mathcal{C}})$. 
In~particular, if the induced cover is principal and $\bigcap_{i=1}^n C_i=\emptyset$ then $\mathcal{C}$ is principal. (Confer~Lemma~\ref{lem:G-KKMmaps} below.)
\end{itemize}
\end{rmrk}

We now have all the definitions needed to replace spheres in Theorem~\ref{thm:KKM} with general `noncontractible' spaces.

\begin{thm}\label{thm:Gg1}
Let $X$ be a compact normal space,
let $M$ be a contractible metric space, 
and let $f\colon X\to M$ be a continuous map. 
Then for any non--null--homotopic cover~$\mathcal C$ of~$X$,
a set of spherical $f$--neighbors intersects at least $\rk(\mathcal C)+1$ elements of~$\mathcal C$.
In~particular, for any principal cover,
a set of spherical $f$--neighbors intersects all elements of the cover.
\end{thm}

Theorem~\ref{thm:Gg1} implies Theorem~\ref{thm:KKM} because 
each KKM cover either is principal or all of its elements have a common point; furthermore, the maps in the homotopy class $[\mathcal{C}]$ corresponding to each principal KKM cover $\mathcal{C}$ are of degree one, so that $[\mathcal{C}]$ contains a homeomorphism (see Corollaries 2.1--2.3 in~\cite{Mus17}).

\begin{rmrk}
Note that $X$ in Theorem~\ref{thm:Gg1} is not assumed to be connected. 
Figure~\ref{fig:touches-plus-point} shows an example with $X=\SSS^1\cup\{b',c'\}$ (cf.~Fig.~\ref{fig:touches}).
\begin{figure}[h]
\scalebox{0.6}{
\begin{tikzpicture}
\draw[thick,rounded corners=20pt,densely dotted]
(0,-3) -- (0,0) -- (1,3)-- (5,4) -- (6,2) -- (1.45,1.35); 
\draw[thick,rounded corners=20pt]
(1.45,1.35) -- (-1,1) -- (-1, -1) -- (7,-2)--(8,-1);
\draw[thick,rounded corners=20pt,densely dashed]
(8,-1) -- (9,0) -- (8,5) -- (6,5) -- (4, -3) -- (0, -4) -- (0,-3);
\draw[thick, densely dashed] (3,0.03) circle (0.1) node [right,black]{\large $f(c')$};
\fill (1,-2.73) circle (0.1) node [right,black]{\large $f(b')$};
\draw[thick,gray] (0.5,-3) circle (0.48) node[black]{\large $D^2$};
\draw (3.16,4.03)node{\large $f(a)$};
\draw (-1.66,0.03)node{\large $f(b)$};
\draw (9.06,3.03)node{\large $f(c)$};
\draw[thick,densely dotted] (11,4)--(12,4) node[right]{\large $a$};
\draw[thick]                (11,3.5)--(12,3.5) node[right]{\large $b, b'$};
\draw[thick,densely dashed] (11,3)--(12,3) node[right]{\large $c, c'$};
\end{tikzpicture}}
\caption{An example with disconnected $X=\SSS^1\cup\{b',c'\}$}
\label{fig:touches-plus-point}
\end{figure}
\end{rmrk}

\begin{rmrk}
Combining the idea that $X$ in Theorem~\ref{thm:Gg1} is not necessarily connected, with switching attention to the image of the cover, leads to generalizations of Helly's theorem and the KKM lemma. See also Lemma~\ref{lem:G-KKMmaps} below. We do not develop this line in the present paper. 
\end{rmrk}

\subsection{EP triples and ranks, and spherical $f$--neighbors}

We are going to upgrade Theorem~\ref{thm:Gg1} to a more general Theorem~\ref{thm:G-EPtriple}, which covers the case of maps to not necessarily contractible spaces. 
In order to state and prove Theorem~\ref{thm:G-EPtriple}, we introduce concepts of Eilenberg--Pontryagin and Knaster--Kuratowski--Mazurkiewicz ranks.

\begin{defn}[Eilenberg--Pontryagin triples and ranks] 
\label{def:G-EP}
Let $Z$ be a topological space with a subspace~$A$, let $K$ be a finite simplicial complex, and let $[h]$ be a homotopy class in~$[A,K]$.
We say that $(Z,A,[h])$ is an \emph{Eilenberg--Pontryagin triple} if no map in~$[h]$ extends to a continuous map $Z\to K$.

We define the \emph{Eilenberg--Pontryagin rank} (\emph{EP rank}) $\rk(Z,A,[h])$ of the triple $(Z,A,[h])$ to be the least nonnegative integer $k$ such that there exists a continuous map $H\colon Z\to K\cup\sk_k(\Delta_K)$ whose restriction $H|_A$ is homotopic in $K\cup\sk_k(\Delta_K)$ to the maps of~$[h]$, where $\Delta_K$ is the simplex spanned by the vertices of~$K$ and containing~$K$ as a subcomplex. 
(Since $\Delta_K$ is contractible, the EP rank is well defined and does not exceed the dimension of~$\Delta_K$.)
\end{defn}

\begin{rmrk}
In terms of Definition~\ref{def:G-EP}, a triple $(Z,A,[h])$ is Eilenberg--Pontryagin if and only if it is of nonzero EP rank (because $K\cup\sk_0(\Delta_K)=K$).
\end{rmrk}

\begin{rmrk}
\label{rem:EP-rk-hmt-rk}
Since any constant map extends to any ambient space, 
it follows that in terms of Definitions~\ref{def:homotopy-ranks} and~\ref{def:G-EP}, for any 
$Z$, $A$, $K$, and $h$ we have 
$$
\rk(Z,A,[h])~\le~\rk(h).
$$
Furthermore, if $A$ is contractible in~$Z$, then we have
$$
\rk(Z,A,[h])~=~\rk(h).
$$
In particular, if $\mathcal{C}$ is a finite closed cover of~$A$ and $A$ is contractible in~$Z$, then
$$
\rk(Z,A,[\mathcal{C}])~=~\rk(\mathcal{C}).
$$
\end{rmrk}

\begin{example} \ 
\begin{itemize}
\item If $Z=\Delta^n$, $A=K=\partial\Delta^n$, and $h=\operatorname{id}$, then $\rk(Z,A,[h])=\rk(h)=n$.
\item If $Z=K=\Delta^n$, $A=\partial\Delta^n$, and $h=\operatorname{id}$, then $\rk(Z,A,[h])=\rk(h)=0$.
\end{itemize}
\end{example}

\begin{example}
We have $\rk(Z,A,[h])=0$ whenever $A$ is a retract of~$Z$.
\end{example}

\begin{example}
Let $W$ be an orientable compact PL $m$-manifold with connected nonempty boundary $\partial W$, and let $h\colon\partial W \to \partial\Delta^n$ be a continuous map. Then we have $\rk(h)\in\{0,n\}$, $\rk(W,\partial W,[h])\in\{0,n\}$, and $\rk(W,\partial W,[h])\le\rk(h)$.
\begin{itemize}
\item If $m=n$, then $\rk(W,\partial W,[h])=\rk(h)$ (this follows from the Hopf degree theorem; see the proof of Corollary~\ref{cor:manifold} below).

\item If $W=\Delta^m$, then $\rk(W,\partial W,[h])=\rk(h)$ (because $\Delta^m$ is contractible; see Remark~\ref{rem:EP-rk-hmt-rk}).

\item Results of \cite{MW18} imply however that for any $m$ and $n$ with nontrivial $\pi_{m-1}(\SSS^{n-1})$ and for any non--null--homotopic 
$h\colon \SSS^{m-1}\to\partial\Delta^n$ there exists an $m$-manifold $W$ with $\partial W=\SSS^{m-1}$ such that $h$ extends to a continuous map $W\to\partial\Delta^n$,
so that we have $\rk(W,\partial W,[h])=0$ and $\rk(h)=n$. 
\end{itemize}
%
\end{example}

\begin{defn}[Knaster--Kuratowski--Mazurkiewicz rank] 
\label{def:G-KKM}
Let $Z$ be a topological space, and let $\mathcal{S}=\{S_1,\ldots,S_n\}$ be a collection of subsets in~$Z$.
We say that the pair $(Z,\mathcal S)$ is a \emph{Knaster--Kuratowski--Mazurkiewicz} (\emph{KKM}) \emph{system} if no closed cover $\{E_1,\ldots,E_n\}$ of~$Z$ with $S_i\subset E_i$ for all~$i$ has the same nerve as~$\mathcal{S}$.

We define the \emph{KKM rank} $\rk(Z,\mathcal S)$ of the pair $(Z,\mathcal S)$ to be the least integer~$k$ such that there exists a closed cover $\mathcal{E}=\{E_1,\ldots,E_n\}$ of~$Z$ with $S_i\subset E_i$ for all~$i$ such that the dimension of $\mathcal{N}(\mathcal{E})\setminus \mathcal{N}(\mathcal{S})$ is~$k$. 
\end{defn}

\begin{rmrk}
In terms of Definition~\ref{def:G-KKM}, a pair $(Z,\mathcal S)$ is a KKM system if and only if it is of nonzero KKM rank.
\end{rmrk}

\begin{example}
If $Z=\Delta^m$ and $\mathcal{S}=\{S_1,\ldots,S_{m+1}\}$ is a KKM cover of~$\partial \Delta^m$, then $\rk(Z,\mathcal S)=m$ by the KKM lemma.
\end{example}

\begin{example}
We have $\rk(Z,\mathcal S)=0$ whenever $\mathcal{S}$ is a closed cover of a retract of~$Z$.
\end{example}

\begin{example}
We have $\rk(Z,\{S_1,\ldots,S_n\})=0$ whenever $\bigcap_{i=1}^n S_i\neq\emptyset$.
\end{example}

\begin{lem}
\label{thm:G-EP-KKM}
Let a normal space $Z$ contain a normal space $A$ as a subspace,
let $\mathcal C$ be a closed cover of~$A$, and let $[\mathcal C]$ be the corresponding homotopy class in~$[A,\mathcal{N}(\mathcal{C})]$, where $\mathcal{N}(\mathcal{C})$ is the nerve. 
Then 
the EP rank of the triple $(Z,A,[\mathcal C])$
does not exceed
the KKM rank of the system $(Z,\mathcal C)$\textup{:}
$$\rk(Z,A,[\mathcal C])\le\rk(Z,\mathcal C).$$
Furthermore, if~$A$ is closed in~$Z$, then 
$$\rk(Z,A,[\mathcal C])=\rk(Z,\mathcal C).$$
\end{lem}

Lemma~\ref{thm:G-EP-KKM} is proved in the next section. 

\begin{example}[Showing that the closedness requirement of $A$ in the second part of Lemma~\ref{thm:G-EP-KKM} is essential]
If $X$ is a compact normal space, $\mathcal C'=\{C'_1,\ldots,C'_n\}$ is a closed cover of~$X$ with $\bigcap_{i=1}^n C_i=\emptyset$ and each $C'_i$ nonempty, $Z$ is the cone over~$X$, $z_0$ is the top of~$Z$, $A=Z\setminus\{z_0\}$, $C''_i$ is the subcone in~$Z$ over $C'_i$, $C_i=C''_i\setminus \{z_0\}$, and $\mathcal C=\{C_1,\ldots,C_n\}$, then the KKM rank $\rk(Z,\mathcal{C})$ is $n-1$ and the EP rank $\rk(Z,A,[\mathcal C])$ is one more than the dimension of the nerve~$\mathcal{N}(\mathcal{C}')$. (See Fig.~\ref{fig:wheel} with $X=\SSS^1$.)
For example, if $n>2$ and the elements of $\mathcal C'$ are pairwise disjoint,
then 
$$\rk(Z,A,[\mathcal C])=1<n-1=\rk(Z,\mathcal C).$$
\begin{figure}[ht]
\begin{tikzpicture}
\draw (0,0) circle (2cm) node{$*$}; 
\draw (3mm,-1mm) node{$z_0$};
\draw (0cm,0cm) -- (90:2cm); 
\draw (0cm,0cm) -- (90-72:2cm); 
\draw (0cm,0cm) -- (90-72-72:2cm); 
\draw (0cm,0cm) -- (90+72:2cm); 
\draw (0cm,0cm) -- (90+72+72:2cm); 
\draw (90-36:12mm) node{$C_1$}; 
\draw (90-36-72:12mm) node{$C_2$}; 
\draw (90-36-72-72:12mm) node{$C_3$}; 
\draw (90-36-72-72-72:12mm) node{$C_4$}; 
\draw (90-36-72-72-72-72:12mm) node{$C_5$}; 
\end{tikzpicture}
\caption{A disk with $\rk(Z,A,[\mathcal C])=2<4=\rk(Z,\mathcal C)$}
\label{fig:wheel}
\end{figure}
\end{example}

\begin{rmrk}
\label{rem:KKM-rk-hmt-rk}
Lemma~\ref{thm:G-EP-KKM} implies (see Remark~\ref{rem:EP-rk-hmt-rk}) that, given a compact normal space $A$ with a finite closed cover~$\mathcal{C}$, for any ambient normal space $Z\supset A$ we have
$
\rk(Z,A,[\mathcal C])=\rk(Z,\mathcal{C})\le\rk(\mathcal{C}),
$
and 
$
\rk(Z,A,[\mathcal C])=\rk(Z,\mathcal{C})=\rk(\mathcal{C})
$
if $A$ is contractible in~$Z$. (This generalizes Theorem~2.2 from~\cite{Mus16}.)
\end{rmrk}

\begin{thm}\label{thm:G-EPtriple}
Let $A$ be a compact normal space, 
let $\mathcal C$ be a closed cover of~$A$,
and let $[\mathcal C]$ be the corresponding homotopy class in~$[A,\mathcal{N}(\mathcal{C})]$, where $\mathcal{N}(\mathcal{C})$ is the nerve. 
Let $Z$ be a normal space containing~$A$ as a subspace.
If the triple $(Z,A,[\mathcal C])$ is Eilenberg--Pontryagin, with EP rank $\rk(Z,A,[\mathcal C])>0$,
then for any metric space~$M$ and any continuous map $f\colon A\to M$ that extends to a continuous map $Z\to M$, a set of spherical $f$--neighbors intersects at least $\rk(Z,A,[\mathcal C])+1$ elements of~$\mathcal C$.
\end{thm}

Theorem~\ref{thm:G-EPtriple} is proved in the next section. 

\begin{proof}[Proof of Theorem~\ref{thm:Gg1}]
We deduce Theorem~\ref{thm:Gg1} from Theorem~\ref{thm:G-EPtriple}.
Let $X$, $M$, $f$, and~$\mathcal C$ be as in Theorem~\ref{thm:Gg1}.
Set $\operatorname{Cone}(X):=(X\times [0,1])/(X\times\{0\})$ and identify $X$ with~$X\times\{1\}\subset \operatorname{Cone}(X)$.
(The cone is normal because $X$ is compact and normal; see, e.\,g., \cite{Nob20}.)
Definitions of ranks imply (see Remark~\ref{rem:EP-rk-hmt-rk}) that
\be
\tag{***}
\label{eq:rank=rank}
\rk(\operatorname{Cone}(X),X,[\mathcal C])=\rk(\mathcal{C}).
\ee
In particular, 
the triple $(\operatorname{Cone}(X),X,[\mathcal C])$ is Eilenberg--Pontryagin since $\mathcal{C}$ is non--null--homotopic. Since $M$ is contractible, it follows that $f$ extends to a continuous map $F\colon \operatorname{Cone}(X)\to M$. We apply Theorem~\ref{thm:G-EPtriple} to the Eilenberg--Pontryagin triple $(\operatorname{Cone}(X),X,[\mathcal C])$, with $Z=\operatorname{Cone}(X)$ and $A=X$ in the notation of Theorem~\ref{thm:G-EPtriple}, and see that a set of spherical $f$--neighbors intersects at least $\rk(\operatorname{Cone}(X),X,[\mathcal C]) + 1$ elements of~$\mathcal C$. Then Theorem~\ref{thm:Gg1} follows by~\eqref{eq:rank=rank}. 
\end{proof}

\begin{rmrk}
Theorem~\ref{thm:G-EPtriple} has further refinements regarding the number of distinct sets of spherical $f$--neighbors intersecting the prescribed number of cover elements, but we do not develop this line in the present paper. 
\end{rmrk}

\subsection{Corollaries}
Next, we list several corollaries of Theorems~\ref{thm:simplex}--\ref{thm:G-EPtriple}.
In fact, all of the following corollaries, except for Corollary~\ref{cor:manifold}, follow from Theorem~\ref{thm:simplex}.

\begin{defn}
A~continuous map $f\colon  A\to Y$ of an orientable connected closed PL manifold~$A$ to a space~$Y$ is said to be \emph{null--cobordant} if there exists an orientable compact PL manifold~$W$ with $\partial W = A$ and a continuous map $F\colon W\to Y$ such that $F|_A=f$.
\end{defn}

\begin{cor}[{cf.~\cite[Theorem~2.6]{Mus16}}]
\label{cor:manifold}
%
Let $A$ be an orientable connected closed PL $n$-manifold, and let $\mathcal C$ be a non--null--homotopic cover of~$A$ such that the nerve of $\mathcal C$ is homeomorphic to the~$n$-sphere.
Then for any metric space~$M$ and any null--cobordant map $f\colon A\to M$, 
a set of spherical $f$--neighbors intersects at least $n+2$ elements of~$\mathcal C$.
In~particular, if $\mathcal C$ is principal and contains precisely $n+2$ elements, then a set of spherical $f$--neighbors intersects all elements of~$\mathcal C$.
\end{cor}

\begin{proof}
If $f\colon A\to M$ is null--cobordant, then there are an orientable compact PL $(n+1)$-manifold~$Z$ with $\partial Z = A$ and a continuous map $F\colon Z\to M$ with $F|_A=f$.
A homological argument shows that for each continuous map $H\colon Z\to\mathcal{N}(\mathcal{C})\cong\SSS^n$, the restriction $H|_A$ is of zero degree. 
Then the Hopf degree theorem implies that $H|_A$ is null--homotopic. 
This means that the triple $(Z,A,[\mathcal C])$ is Eilenberg--Pontryagin and the statement follows by Theorem~\ref{thm:G-EPtriple}.
\end{proof}

\begin{rmrk}[The dimensional restriction in Corollary~\ref{cor:manifold} is essential]
It is shown in~\cite{MW18} that any continuous map $\SSS^m\to\SSS^n$ is null--cobordant if $m>n$.
Let $m$ and $n$ be such that $m>n$ and $\pi_m(\SSS^n)$ is nontrivial, and let $h\colon\SSS^m\to\partial\Delta^{n+1}$ be a non--null--homotopic continuous map. Then there exists an orientable compact PL $(m+1)$-manifold $W$ with $\partial W=\SSS^m$ and a continuous map $H\colon W\to\partial\Delta^{n+1}$ such that $H|_{\partial W}=h$. Let $\mathcal{C}$ be the closed cover of~$\partial W$ composed of the inverse images of the facets of~$\Delta^{n+1}$. Then $[\mathcal{C}]=[h]$ and~$\mathcal{C}$ is principal.
We embed $W$ into a Euclidean ball $B^N$ of large dimension and `tiny' diameter, then embed $W$ into the product $\partial\Delta^{n+1}\times B^N$ such that the projection of this embedding to $\partial\Delta^{n+1}$ yields $H$, and take the induced metric on~$W$. 
Now, let $f\colon \partial W\to W$ be the identity map. Then $f$ is null--cobordant but no
set of spherical $f$--neighbors intersects all elements of~$\mathcal C$ if the diameter of $B^N$ is sufficiently small. 
\end{rmrk}



\begin{cor}
\label{thm2}
Let $M$ be a contractible metric space, let $\SSS^n$ be the Euclidean unit $n$\nobreakdash-sphere  
in Euclidean $(n+1)$-space~$\R^{n+1}$, and let $f\colon\SSS^n\to M$ be a continuous map. 
Then there exists a pair $\{p,q\}$ of spherical $f$--neighbors 
such that 
the Euclidean distance $\|p-q\|$ 
is at least~$\sqrt{\frac{n+2}{n}}$.
\end{cor}

Corollary~\ref{thm2} is proved in the next section. 

\begin{rmrk}
In~\cite{MM20}
we show that if~$M=\R^{m}$ with $m>n$, then the constant~$\sqrt{{(n+2)}/{n}}$ in Corollary~\ref{thm2} 
(the Euclidean distance between the centers of adjacent $(n-1)$-simplices of the regular triangulation of~$\SSS^n$)
can be replaced with 
$\sqrt{{2(n+2)}/(n+1)}$ (the Euclidean distance between vertices of the regular triangulation of~$\SSS^n$),   
which is the best possible.  
Our proof for the Euclidean case $M=\R^{m}$ is based on the Delaunay triangulations and we do not know whether it extends to all contractible~$M$.
\end{rmrk}

\begin{cor}
\label{cor:with-nerve-theorem}
Let $M$ be a contractible metric space, 
let $P$ be a convex $n$-polytope, 
and let $f\colon\dd P \to M$ be a continuous map. 
Then a set of spherical $f$--neighbors intersects at least $n+1$ facets of~$P$.
\end{cor}

\begin{proof}[Proof via Theorem~\ref{thm:simplex}]
Corollary~\ref{cor:with-nerve-theorem} is an `equivalent generalization' of Theorem~\ref{thm:simplex} because the $(n-2)$-skeleton of any convex $n$-polytope contains the $(n-2)$-skeleton of the $n$-simplex as a topological subspace (see~\cite{Gru65}).
\end{proof}

\begin{proof}[Proof via Theorem~\ref{thm:Gg1}]
Clearly, the cover~$\mathcal{C}$ of~$\partial P$ composed of the facets of~$P$ is non--null--homotopic of rank~$n$ because $\mathcal{C}$ is a \emph{good} cover (that is, any intersection of elements in $\mathcal{C}$ is contractible), so the nerve of~$\mathcal{C}$ has homotopy type of~$\partial P\cong\SSS^{n-1}$ by the nerve theorem, while the maps in the class~$[\mathcal{C}]$ are homotopy equivalences.
Then a set of spherical $f$--neighbors intersects at least $n+1$ facets of~$P$ by Theorem~\ref{thm:Gg1}.
\end{proof}

Since any collection of $n+1$ facets of the $n$-cube contains a pair of antipodal facets, Corollary~\ref{cor:with-nerve-theorem} implies the following.

\begin{cor}
\label{cor:cube}
Let $M$ be a contractible metric space, 
let $\dd [0,1]^m$ be the boundary of the $m$-dimensional cube $[0,1]^m$, 
and let $f\colon\dd [0,1]^m \to M$ be a continuous map. 
Then there is a pair of spherical $f$--neighbors intersecting antipodal facets of $[0,1]^m$.
\end{cor}

There exists an example of continuous map~$\SSS^2\to\R^3$ showing that 
the statement of Corollary~\ref{cor:cube} about spherical $f$--neighbors lying on antipodal facets holds for neither regular octahedron nor regular dodecahedron nor regular icosahedron.  
A~weaker version of Corollary~\ref{cor:cube} where `antipodal' is replaced with `disjoint' holds for many polytopes.  

\subsection{Radon type theorems}

\begin{defn}[Weak Radon polytopes] 
\label{def:wRp}
We say that 
an $n$-polytope~$P$ 
is \emph{weakly Radon} if for any continuous map $f\colon \dd P\to M$ into any contractible metric space~$M$ there is a pair of spherical $f$--neighbors intersecting two disjoint faces of~$P$.
\end{defn}

We recall some standard definitions.
A~\emph{flag} polytope is a convex polytope such that every collection of pairwise intersecting facets has a nonempty intersection.
A (\emph{combinatorial}) \emph{fullerene} is a simple $3$-polytope with all facets pentagons and hexagons.

A `visual' simply checked sufficient condition for weakly Radon polytopes is provided by the so-called \emph{belts}.
A~\emph{$k$-belt} (or a \emph{prismatic $k$-circuit}) in a $3$-polytope is a cyclic sequence $(F_1, \dots, F_k)$ of $k\ge 3$ facets in which pairs of consecutive facets (including $F_k, F_1$) are adjacent, other pairs of facets do not intersect, and no three facets have a common vertex. 

\begin{cor}
\label{cor:icos2}
\label{cor:belt}
\ 
\begin{enumerate}
\item If the $(n-2)$-skeleton of a convex $n$-polytope~$P$ contains the $(n-2)$-skeleton of the $n$-cube as a topological subspace, then $P$ is weakly Radon. 
\item Each convex $3$-polytope having a $k$-belt with $k\ge4$ is weakly Radon. 
\item Each flag $3$-polytope is weakly Radon. 
\item Each fullerene is weakly Radon. 
\item The regular dodecahedron and the regular icosahedron are weakly Radon. 
\end{enumerate}
\end{cor}

\begin{proof}
Assertion~(1)
follows from Corollary~\ref{cor:cube} in an obvious way. Assertions~(2) and~(5) directly follow from assertion~(1). 
Assertion~(3) follows from assertion~(2) of Corollary~\ref{cor:fullerene} given below.
Assertion~(4) is a particular case of assertion~(3).
\end{proof}

\begin{defn}[Weak Radon rank] 
\label{def:wRr}
If $P$ is an $n$-polytope, $Y$ is a metric space, and $f\colon \dd P\to M$ is a map, we say that two facets~$F_1$ and $F_2$ of~$P$ are \emph{spherical $f$--neighbors} (or that the pair $\{F_1,F_2\}$ is a \emph{pair of spherical $f$--neighbors}) if there is a pair $\{p,q\}$ of spherical $f$--neighbors with $p\in F_1$ and $q\in F_2$.
We say that $f\colon \dd P\to M$ \emph{has weak Radon rank $m$} if 
there are exactly $m$ distinct pairs of facets of~$P$ such that each of these pairs is a pair of disjoint spherical $f$--neighbors.
By the \emph{weak Radon rank} of a polytope~$P$ we mean the least of the weak Radon ranks of continuous maps $f\colon \dd P\to M$ into contractible metric spaces.
\end{defn}

Corollary~\ref{cor:cube} allows us to obtain rough lower bounds on the weak Radon rank.

\begin{defn}[Cubic hemisphere] 
Let $H$ be a subset of the boundary~$\dd P$ of a convex $n$-polytope~$P$. We say that $H$ is a \emph{cubic hemisphere} if there exists a homeomorphism $h\colon [0,1]^n\to P$ such that the restriction of~$h$ to the $(n-2)$-skeleton of~$[0,1]^n$ is a topological embedding to the $(n-2)$-skeleton of~$P$ and $H$ is the image of the union of $n$ facets of~$[0,1]^n$ that have a common vertex.  
\end{defn}

\begin{defn}[Lighthouse independence number] 
\label{def:lin}
We say that a set $Z$ of vertices of an $n$-polytope is \emph{lighthouse independent} if no two vertices in~$Z$ share a facet (equivalently, the corresponding facets of the dual polytope are pairwise disjoint). The \emph{lighthouse independence} number of an $n$-polytope~$P$, $\operatorname{lin}(P)$, is the cardinality of a largest lighthouse independent set of~$P$. 
\end{defn}

\begin{rmrk}
The lighthouse independence number of an $n$-polytope equals the cardinality of a largest set of pairwise disjoint facets of the dual polytope.
\end{rmrk}

\begin{cor}
\label{cor:hemi} 
\label{cor:fullerene}
\
\begin{enumerate}
\item
Let $P$ be a convex $n$-polytope.
If $\dd P$ contains $k$ cubic hemispheres with pairwise disjoint interiors, then the weak Radon rank of~$P$ is at least~$k/2$.
\item
Let $P$ be a flag $3$-polytope {\textup(}e.\,g., a fullerene{\textup)}. Then the weak Radon rank of~$P$ is at least
half the lighthouse independence number of~$P$.
\item 
If $P$ is a flag simple $3$-polytope with $\psi$ facets and $g$ is the largest number of edges in a facet of~$P$, then the 
weak Radon rank of~$P$ is at least
$$
\left.\left\lfloor\frac{2\psi-7}{3g-8}\right\rfloor\right/2,
$$
where $\lfloor \cdot\rfloor$ stands for the floor function.
\item If $P$ is a fullerene with $\psi$ facets, then the weak Radon rank of~$P$ is at least
$$
\left.\left\lfloor\frac{\psi-3}{5}\right\rfloor\right/2.
$$
\item The weak Radon rank of the regular dodecahedron is at least~$2$.
\item The weak Radon rank of the regular icosahedron is at least~$2$.
\item The weak Radon rank of the cube is~$1$.
\end{enumerate}
\end{cor}

Corollary~\ref{cor:hemi} is proved in the next section. 

\section{Proofs}

\begin{proof} [Proof of Lemma~\ref{thm:G-EP-KKM}]
1) We show that $\rk(Z,A,[\mathcal C])\le \rk(Z,\mathcal C)$.

Let $\mathcal C=\{C_1,\ldots,C_n\}$, let $\Delta^{n-1}$ denote the simplex spanned by the vertices of the nerve~$\mathcal{N}(\mathcal{C})$ of~$\mathcal C$ so that $\mathcal{N}(\mathcal{C})$ is a subcomplex of~$\Delta^{n-1}$, and let $r:=\rk(Z,\mathcal C)$.
By the definition of the KKM rank there exists a closed cover $\mathcal{E}=\{E_1,\ldots,E_n\}$ of~$Z$ with $C_i\subset E_i$ for all~$i$ such that the dimension of $\mathcal{N}(\mathcal{E})\setminus \mathcal{N}(\mathcal{C})$ is~$r$. 
Therefore, the union $\mathcal{N}(\mathcal{C})\cup\sk_r(\Delta^{n-1})$ contains $\mathcal{N}(\mathcal{E})$.
Set
$$
\mathcal{E}_A:=\{E_1\cap A, \dots, E_n\cap A\}.
$$ 
Since $C_i\subset E_i$ for all~$i$, it follows that the nerve $\mathcal{N}(\mathcal{E}_A)$ contains
$\mathcal{N}(\mathcal{C})$.
We have
$$
\mathcal{N}(\mathcal{C})~\subset~ \mathcal{N}(\mathcal{E}_A) ~\subset~ \mathcal{N}(\mathcal{E}) ~\subset~ \mathcal{N}(\mathcal{C})\cup\sk_r(\Delta^{n-1}).
$$
Let $[\mathcal{E}]$ be the homotopy class in~$[Z,\mathcal{N}(\mathcal{E})]$ determined by~$\mathcal E$,
and let $F\colon Z\to \mathcal{N}(\mathcal{E})$ be a map in $[\mathcal{E}]$.
Let $[\mathcal{E}_A]$ be the homotopy class in~$[A,\mathcal{N}(\mathcal{E}_A)]$ determined by~$\mathcal{E}_A$, and let $f'\colon A \to \mathcal{N}(\mathcal{E}_A)$ be a map in $[\mathcal{E}_A]$.
Let $f\colon A \to \mathcal{N}(\mathcal{C})$ be a map in $[\mathcal{C}]\in[A,\mathcal{N}(\mathcal{C})]$.

Since $E_i\cap A$ contains~$C_i$ for each $i$, the argument preceding Definition~\ref{def2} shows that $f$ and $f'$ are homotopic in~$\mathcal{N}(\mathcal{E}_A)$.
By construction, $F|_A$ and $f'$ are homotopic in~$\mathcal{N}(\mathcal{E})$.
Thus, $F|_A$ and $f$ are homotopic in~$\mathcal{N}(\mathcal{E})$ and hence in~$\mathcal{N}(\mathcal{C})\cup\sk_r(\Delta^{n-1})$ as well.
By the definition of the EP rank this means that $\rk(Z,A,[\mathcal C])\le r = \rk(Z,\mathcal C)$.

2) We show that $\rk(Z,\mathcal C)\le \rk(Z,A,[\mathcal C])$ whenever~$A$ is closed in~$Z$.

We start by constructing a specific map $A\to\mathcal{N}(\mathcal{C})$ from the class~$[\mathcal C]$.
Let $\mathcal C=\{C_1,\ldots,C_n\}$, and let $\mathcal{N}(\mathcal{C})$ be a subcomplex in~$\Delta^{n-1}$ (as in the first part of the proof).
Since $A$ is normal, there exists an open cover $\mathcal U=\{U_1,\ldots,U_n\}$ of~$A$ such that $U_i$ contains~$C_i$ for each~$i$ and the nerve of~$\mathcal U$ coincides with that of~$\mathcal C$ (see, e.\,g.,~\cite[Theorem~1.3]{Mor50} and \cite[pp.~31--33]{Iva20}).
The Urysohn lemma for normal spaces implies that for each~$i$ there exists a continuous function $f_i\colon A\to [0,1]$ with $f_i(C_i)=1$ and $f_i(A\setminus U_i)=0$.
Then $\Phi=\{\varphi_1,\ldots,\varphi_n\}$, where $\varphi_i:=f_i/\sum_j f_j$, is a partition of unity subordinate to $\mathcal U$ and such that $\varphi_i^{-1}[1/n,1]$ contains~$C_i$ for all~$i$.
Let
$$
h_{\mathcal U,\Phi}(x):=\sum\limits_{i=1}^n{\varphi_i(x)v_i}
$$
be the corresponding map $A\to\mathcal{N}(\mathcal{C})$ representing the class $[\mathcal C]=[\mathcal U]$ (see the construction preceding Definition~\ref{def2}).

Now, let $p:=\rk(Z,A,[\mathcal C])$.
Then by the definition of the EP rank there exists a continuous map $F\colon Z\to \mathcal{N}(\mathcal{C})\cup\sk_p(\Delta^{n-1})$ such that the restriction $F|_A$ is homotopic to $h_{\mathcal U,\Phi}$ in~$\mathcal{N}(\mathcal{C})\cup\sk_p(\Delta^{n-1})$.
The generalizations of Borsuk's homotopy extension theorem obtained in \cite{Mor75} and \cite{Sta75} imply that, since $A$ is closed in~$Z$, there exists a continuous map $G\colon Z\to \mathcal{N}(\mathcal{C})\cup\sk_p(\Delta^{n-1})$ with $G|_A=h_{\mathcal U,\Phi}$.
Then the collection of subsets $$\mathcal{G}:=\{G_1^{-1}[1/n,1],\ldots,G_n^{-1}[1/n,1]\},$$ where $G_1$, $\ldots$, $G_n$ are the coordinate functions of~$G$, is a closed cover of~$Z$ such that $G_i^{-1}[1/n,1]$ contains $C_i$ for all~$i$ and the nerve $\mathcal{N}(\mathcal{G})$ is contained in $\mathcal{N}(\mathcal{C})\cup\sk_p(\Delta^{n-1})$ so that the dimension of $\mathcal{N}(\mathcal{G})\setminus \mathcal{N}(\mathcal{C})$ is at most~$p$.
By the definition of the KKM rank this means that $\rk(Z,[\mathcal C])\le p = \rk(Z,A,[\mathcal C])$.
\end{proof}

Now we state and prove Lemmas~\ref{lem:G-KKMmaps} and~\ref{lem:G-KKM-metric}, and then deduce Theorem~\ref{thm:G-EPtriple} from Lemmas~\ref{thm:G-EP-KKM}, \ref{lem:G-KKMmaps}, and~\ref{lem:G-KKM-metric}.

\begin{lem}
\label{lem:G-KKMmaps}
Let $(Z,\mathcal C=\{C_1,\ldots,C_n\})$ be a KKM system of rank~$r>0$, 
let $f\colon Z\to Z'$ be a continuous map to a topological space~$Z'$,
and let $\mathcal C'=\{C'_1,\ldots,C'_n\}$ be a family of subsets in~$Z'$ such that $f(C_i)\subset C'_i$ for all~$i$. 
Then either $\{1,\ldots,n\}$ contains a subset $J$ of cardinality $r+1$ such that $\bigcap_{j\in J}C_j=\emptyset$ and $\bigcap_{j\in J}C'_j\neq\emptyset$ or $(Z',\mathcal C')$ is a KKM system of rank at least~$r$.
\end{lem}

\begin{rmrk}
In Lemma~\ref{lem:G-KKMmaps} two key special cases are $C'_i=f(C_i)$ and $f=\operatorname{id}$.
\end{rmrk}

\begin{proof}
If neither $\rk(Z',\mathcal C')\ge r$ nor $\{1,\ldots,n\}$ contains $J$ with $|J|=r+1$ such that $\bigcap_{j\in J}C_j=\emptyset$ and $\bigcap_{j\in J}C'_j\neq\emptyset$, 
then (i) there exists a closed cover $\mathcal{E}'=\{E'_1,\ldots,E'_n\}$ of~$Z'$ with $C'_i\subset E'_i$ for all~$i$ such that the dimension of $\mathcal{N}(\mathcal{E}')\setminus \mathcal{N}(\mathcal{C}')$ is less than~$r$ (by definition) and (ii) the dimension of $\mathcal{N}(\mathcal{C}')\setminus \mathcal{N}(\mathcal{C})$ is less than~$r$.
Consequently, the dimension of $\mathcal{N}(\mathcal{E}')\setminus \mathcal{N}(\mathcal{C})$ is less than~$r$.
The collection $\mathcal{E}=\{E_1,\ldots,E_n\}$ with $E_i:=f^{-1}(E'_i)$ is a closed cover of~$Z$ such that
$C_i\subset E_i$ for all~$i$. The nerve $\mathcal{N}(\mathcal{E}')$ contains $\mathcal{N}(\mathcal{E})$. 
Therefore, the dimension of $\mathcal{N}(\mathcal{E})\setminus \mathcal{N}(\mathcal{C})$ is less than~$r$.
This contradicts the assumption that $r=\rk(Z,\mathcal C)$.
\end{proof}

\begin{lem}
\label{lem:G-KKM-metric}
Let $(Z,\mathcal C=\{C_1,\ldots,C_n\})$ be a KKM system of rank~$r>0$ with metrizable~$Z$ and all $C_i$ compact, and let $d$ be a metric on~$Z$.
Then there exists a closed metric ball 
$$
B_R(x):=\{z\in Z\mid d(z,x)\le R\}, x\in Z, R> 0,
$$
whose interior intersects no element of~$\mathcal{C}$ and whose boundary sphere touches at least $r+1$ elements of~$\mathcal{C}$.
\end{lem}

\begin{proof}
If $z$ is a point and $N$ is a subset in~$Z$, we write
$$
d(z,N):=\inf_{p\in N}d(z,p).
$$
Let $A$ denote the union $\bigcup_{i=1}^{n}C_i$.
For each $i\in\{1,\dots,n\}$, we set
$$
E_i:=\{z\in Z \mid d(z,C_i)=d(z,A)\}.
$$ 
Observe that $E_i$ contains~$C_i$ and is closed so that 
$\{E_1,\dots,E_n\}$ is a closed cover of~$Z$.
Since $(Z,\mathcal C)$ is a KKM system of rank~$r>0$, 
it follows by the definition of KKM rank that 
the set $\{1,\ldots,n\}$ contains a subset $J$ of cardinality $r+1$ such that $\bigcap_{j\in J}C_j=\emptyset$ and $\bigcap_{j\in J}E_j\neq\emptyset$.
Let $x$ be a point in $\bigcap_{j\in J}E_j\neq\emptyset$.
Then the ball $B_{d(x,A)}(x)$ of radius~$d(x,A)$ centered at~$x$
meets the requirements of the lemma (since each of~$C_i$ is compact).
\end{proof}

\begin{proof}[Proof of Theorem~\ref{thm:G-EPtriple}]

Since $(Z,A,[\mathcal C])$ is an Eilenberg--Pontryagin triple, it follows by Lemma~\ref{thm:G-EP-KKM} that $(Z,\mathcal C)$ is a KKM system of rank~$\rk(Z,\mathcal C)=\rk(Z,A,[\mathcal C])$.

Let $\mathcal C=\{C_1,\ldots,C_n\}$. Set $\mathcal{F}:=\{F(C_1),\ldots,F(C_n)\}$.
Then Lemma~\ref{lem:G-KKMmaps} implies that we have two possibilities:
\begin{itemize}
\item[(1)] The dimension of $\mathcal{N}(\mathcal{F})\setminus\mathcal{N}(\mathcal{C})$ is at least~$\rk(Z,\mathcal C)$ so that $\{1,\ldots,n\}$ contains a subset $J$ of cardinality $\rk(Z,\mathcal C)+1$ such that $\bigcap_{j\in J}C_j=\emptyset$ and $\bigcap_{j\in J}F(C_j)\neq\emptyset$. 
\item[(2)] The pair $(M,\mathcal{F})$ is a KKM system of rank at least~$\rk(Z,\mathcal C)$. 
\end{itemize} 

In case (1), for any point $x\in\bigcap_{j\in J}F(C_j)$ the set $F^{-1}(x)$ is a set of spherical $F|_A$--neighbors that intersects all elements of~$\{C_j\}_{j\in J}$, which proves the theorem.

In case (2), the required statement follows by Lemma~\ref{lem:G-KKM-metric} applied to $(M,\mathcal{F})$. 
\end{proof}

\begin{proof}[Proof of Corollary~\ref{thm2}] 
We use a spherical version of Theorem~\ref{thm:simplex}. 
Let $T$ be a regular triangulation of the unit sphere $\SSS^n$, and let $\tilde\Delta_1$, $\ldots$, $\tilde\Delta_{n+2}$, be the $n$-simplices of~$T$: all of $\tilde\Delta_i$ are regular spherical simplices with Euclidean distances between vertices 
\begin{equation}
\label{eq:EEL}
d_{n,Eu}=\sqrt{\frac{2(n+2)}{n+1}}
\end{equation}
and angular edge length 
\begin{equation}
\label{eq:AEL}
d_{n,A}=2\arcsin\frac{d_{n,Eu}}{2}=2\arcsin\sqrt{\frac{n+2}{2(n+1)}}=\arccos\left(\frac{-1}{n+1}\right).
\end{equation} 

We recall that the \emph{circumradius} of a compact set~$Q$ in a metric space is defined as the radius of a least metric ball containing~$Q$. 
If $Q$ is a compact subset of~$\SSS^n$
we denote by $\crA Q$ and $\diamA Q$, respectively, the circumradius and diameter of~$Q$ with respect to the angular metric, and $\diamE Q$ will stand for the Euclidean diameter of~$Q$ in~$\R^{n+1}\supset\SSS^n$.
Under this notation Dekster's extension~\cite{Dekster} of the Jung theorem says that for any compact subset $Q$ of~$\SSS^n$ we have
\begin{equation*}
\label{eq:Dekster}
2\arcsin\left(\sqrt{\frac{n+1}{2n}} \sin (\crA Q)\right)~\le~\diamA Q.
\end{equation*}
This immediately implies that in the case where $\crA Q\le \pi/2$ we have
\begin{equation}
\label{eq:DeksterE}
\sqrt{\frac{2(n+1)}{n}} \sin (\crA Q)~\le~\diamE Q.
\end{equation}

Another auxiliary fact we need is that 
\begin{equation}
\label{eqn:facet-diameter}
\diamA \tilde\Delta_i = \pi-d_{n,A}/2.
\end{equation}
Indeed, observe that $\tilde\Delta_i$ is the intersection of a finite number of closed hemispheres and hence its boundary is composed of fragments of great hyperspheres, which are geodesic in~$\SSS^n$.
Therefore, if $a$ and $b$ are two points in $\tilde\Delta_i$ such that neither $a$ nor $b$ is a vertex of $\tilde\Delta_i$, then $\diamA\{a,b\}<\diamA \tilde\Delta_i$ because $\tilde\Delta_i$ contains two geodesic arcs\footnote{By geodesic arcs in~$\SSS^n$ we mean arcs of great circles.} $\alpha$ and $\beta$ such that 
$\alpha$ contains $a$ in its relative interior and $\beta$ contains $b$ in its relative interior. Since $\tilde\Delta_i$ is contained in the interior of a hemisphere so that $a$ and $b$ are not antipodal, it follows that there exist $a'\in\alpha\subset\tilde\Delta_i$ and $b'\in\beta\subset\tilde\Delta_i$ with $\diamA\{a,b\}<\diamA \{a',b'\}$ (imagine the interposition of $\alpha$, $\beta$, and the metric ball $D\subset\SSS^n$ of diameter $\diamA\{a,b\}$ containing $a$ and~$b$). Thus, if $a$ and $b$ are two points in $\tilde\Delta_i$ such that $\diamA\{a,b\}=\diamA \tilde\Delta_i$, then one of $a$ and $b$ is a vertex of $\tilde\Delta_i$ and we easily obtain~\eqref{eqn:facet-diameter} by considering the regular triangulation of $\SSS^n$ dual (antipodal) to~$T$.

Now, we pass to the proof of Corollary~\ref{thm2}.
If we have a continuous map $f\colon \SSS^n\to M$, 
then Theorem~\ref{thm:simplex} implies that a finite set $\mathscr{P}$ of spherical $f$--neighbors intersects all of~$\tilde\Delta_i$.
We need to prove that 
\begin{equation}
\label{eqn:required}
\diamE \mathscr{P}~\ge~\sqrt{\frac{n+2}{n}}.
\end{equation}

Let~$B\subset\SSS^n$ be a metric ball with angular radius $\crA \mathscr{P}$ containing~$\mathscr{P}$, 
let $C\in\SSS^n$ be the center of~$B$, 
let $A\in\SSS^n$ be the antipode of~$C$, 
let $\tilde\Delta_k$ be a simplex of~$T$ containing~$A$, 
and let $B_2\subset \SSS^n$ be the metric ball centred at~$A$ of angular radius~$\diamA \tilde\Delta_k = \pi-d_{n,A}/2$ (see \eqref{eqn:facet-diameter}). 
Then $B_2$ contains $\tilde\Delta_k$. 
Since $\mathscr{P}$ intersects $\tilde\Delta_k$ while $\mathscr{P}\subset B$ and $\tilde\Delta_k \subset B_2$, it follows that $B$ intersects~$B_2$. Therefore, we have
\begin{equation}
\label{eqn:key}
\crA \mathscr{P}
~=~
\crA B 
~\ge~ 
\pi- \crA B_2
~=~
d_{n,A}/2.
\end{equation}
The situation splits in two cases:
\begin{itemize}
\item[(i)] $\crA \mathscr{P} > \pi/2$ (i.\,e., no hemisphere contains $\mathscr{P}$);
\item[(ii)] $\crA \mathscr{P} \le \pi/2$.
\end{itemize}

In case (i) we observe that since no hemisphere contains~$\mathscr{P}$, it follows that no Euclidean ball in~$\R^n$ of radius less than $1$ contains~$\mathscr{P}$. 
Then the Jung theorem\footnote{A discussion and materials concerning the Jung theorem and containment in hemispheres see in \cite{Ver52}, \cite{Gus}, \cite[pp.~112, 113, 132--136]{DGK63}, \cite{Tim65}, \cite{Kle}, \cite{Ale77}, and \cite[Proposition~2.4]{AMS07}.}
says that $\diamE \mathscr{P} \ge d_{n,Eu}$, which implies the required~\eqref{eqn:required}.
 
In case (ii), \eqref{eq:DeksterE} is applicable and yields
\begin{multline*}
\diamE \mathscr{P}
~\stackrel{\eqref{eq:DeksterE}}{\ge}~
\sqrt{\frac{2(n+1)}{n}} \cdot \sin (\crA \mathscr{P})
~\stackrel{\eqref{eqn:key},\text{(ii)}}{\ge}~
\sqrt{\frac{2(n+1)}{n}} \cdot \sin \left(\frac{d_{n,A}}{2}\right) \\
~\stackrel{\eqref{eq:AEL}}{=}~
\sqrt{\frac{2(n+1)}{n}} \cdot \frac{d_{n,Eu}}{2}  
~\stackrel{\eqref{eq:EEL}}{=}~
\sqrt{\frac{2(n+1)}{n}} \cdot \sqrt{\frac{n+2}{2(n+1)}} 
~=~
\sqrt{\frac{n+2}{n}}.
\qedhere
\end{multline*}
\end{proof}

\begin{rmrk}
It would be interesting to find a way to upgrade the above proof of Corollary~\ref{thm2} by considering the family of all regular triangulations of the unit sphere~$\SSS^n$.  
\end{rmrk}

\begin{proof}[Proof of Corollary~\ref{cor:fullerene}]
(1) Corollary~\ref{cor:cube} implies that if 
$M$ is a contractible metric space and $f\colon\dd P \to M$ is a continuous map, then each cubic hemisphere in~$\dd P$ contains a facet that is a member of a pair of disjoint facets that are spherical $f$--neighbors.
The statement follows.

(2) Proposition~\ref{prop:cubical} below implies that if the lighthouse independence number of a flag $3$-polytope~$P$ is~$k$, then 
$\dd P$ contains $k$ cubic hemispheres with pairwise disjoint interiors. 
This implies the required assertion by assertion~(1) of the corollary.

Assertions~(3) and~(4) of Corollary~\ref{cor:fullerene} follow from assertion~(2) and Proposition~\ref{thm:lin} below.

Assertions~(5) and~(6) follow from assertion~(2) because the direct check shows that the lighthouse independence number of the regular dodecahedron is~$4$, and the lighthouse independence number of the regular icosahedron is~$3$.

As for the weak Radon rank of the cube (assertion~(7)), Corollary~\ref{cor:cube} shows that it is at least~$1$ and an example where $\dd [0,1]^3$ is mapped to an oblate spheroid in~$\R^3$ shows that it is at most~$1$.
\end{proof}

We say that a vertex~$v$ of a polytope is \emph{cubical} if the union of the facets containing $v$ is a cubic hemisphere.

\begin{prop}
\label{prop:cubical} 
All vertices of a flag $3$-polytope are cubical.
\end{prop}

\begin{proof}[Proof of Proposition~\ref{prop:cubical}]
Let $v$ be a vertex of a flag $3$-polytope~$P$.
Observe that no facet of~$P$ is a triangle (because any triangular facet together with the three adjacent ones form a collection of four pairwise intersecting facets with no common point).
Therefore, each facet of~$P$ containing~$v$ has a vertex that is not adjacent to~$v$.
Let $v_1$, $v_2$, and $v_3$ be three such vertices lying on three distinct facets containing~$v$. Let $D$ denote the union of the facets of~$P$ that do not contain~$v$. Then $D$ is a topological disk with the points~$v_1$, $v_2$, and $v_3$ on its boundary. 
Since $P$ if flag, we see that 
\begin{itemize}
\item no facet contained in~$D$ intersects three of the facets not contained in~$D$, 
\item no facet of~$P$ splits~$D$ (in the sense that $D\setminus F$ is connected for each facet~$F$).
\end{itemize}

This implies that 
\begin{itemize}
\item each of the vertices~$v_1$, $v_2$, and $v_3$ is incident to an edge of~$P$ whose second endpoint is contained in the interior of~$D$ (in particular, the interior of $D$ contains at least one vertex of~$P$), and
\item the subgraph~$G_D$ in the $1$-skeleton~$P_1$ of~$P$ induced by the vertices of~$P$ contained in the interior of~$D$ is connected.
\end{itemize}

Thus, each of $v_1$, $v_2$, and $v_3$ is adjacent to a vertex of the connected subgraph~$G_D$ in~$P_1$.
This easily implies that $P_1$ contains a $Y$-homeomorphic subgraph~$Y'$ 
that is contained in~$D$ and intersects the boundary~$\dd D$ exactly in the set~$\{v_1,v_2,v_3\}$.

Besides, since~$v_1$, $v_2$, and $v_3$ belong to three distinct facets containing~$v$, it follows that there exists a triple of edges in~$P_1$ incident to~$v$ whose endpoints split $\dd D$ into three arcs each of which contains exactly one of~$v_1$, $v_2$, and $v_3$. Clearly, the union of these edges with $\dd D$ and $Y'$ is a graph homeomorphic to the cube $1$-skeleton. 
This shows that $v$ is cubical.
\end{proof}

\begin{prop}
\label{thm:lin} 
\
\begin{enumerate}
\item 
Let $P$ be a flag simple $3$-polytope with $\psi$ facets, and let $g$ be the largest number of edges in a facet of~$P$. Then the 
lighthouse independence number of~$P$ is at least
$$
\left\lfloor\frac{2\psi-7}{3g-8}\right\rfloor.
$$
\item If $P$ is a fullerene with $\psi$ facets, then the lighthouse independence number of~$P$ is at least
$$
\left\lfloor\frac{\psi-3}{5}\right\rfloor.
$$
\end{enumerate}
\end{prop}

\begin{proof}[Proof of Proposition~\ref{thm:lin}]
In the proof if $v$ is a vertex of~$P$, we denote by $L(v)$ the union of facets of~$P$ that contain~$v$.

We construct a lighthouse independent set by the following algorithm.
First we choose a vertex~$v_1$ of $P$ such that the number of vertices in $L(v_1)$ is the least possible and set $W_1=L(v_1)$. The number of vertices in $L(v_1)$ is at most $3g-5$.

At each next step, being given $W_i\subset P$ such that a vertex of~$P$ is not in~$W_i$, we take a vertex~$v_{i+1}$ of $P$ in $P\setminus W_i$ such that the number of vertices in $L(v_{i+1})\setminus W_i$ is the least possible and set $W_{i+1}=W_i\cup L(v_{i+1})$. 
Observe that if a vertex $v$ of~$P$ is not in~$W_i$ and adjacent to a vertex in~$W_i$, then $L(v)$ shares at least $3$ vertices with~$W_i$. This implies that the number of vertices in $L(v_{i+1})\setminus W_i$ is at most $3g-8$.

Therefore, if~$P$ has~$N$ vertices this algorithm produces a lighthouse independent set $v_1$, $v_2$, ... with at least
$$
1+\left\lfloor\frac{N-(3g-5)}{3g-8}\right\rfloor=\left\lfloor\frac{N-3}{3g-8}\right\rfloor
$$
elements.
Since $P$ is simple, Euler's formula yields $N=2\psi-4$. This proves the required estimate.

The case of fullerenes follows if we observe that when $v_1$ is a vertex of a pentagon, then 
the number of vertices in $L(v_1)$ is at most~$12$.
\end{proof}

\section{Concluding remarks}

Now we discuss several concepts and open questions.


\begin{enumerate}

\item 
{\bf The Hopf theorem.} The trefoil curve in Fig.~\ref{fig:trefoil} shows that there exists a continuous map $f\colon \SSS^1\to\R^2$ with no pair of spherical $f$--neighbors having distance less than $\sqrt{3}$ between them.
\begin{figure}[ht]
\scalebox{0.6}{
\begin{tikzpicture}
\draw[thick, yshift=0.577cm] (-1,0) .. controls (-1,0.555) and (-0.555,1) .. (0,1)
               .. controls (0.555,1) and (1,0.555) .. (1,0)
               .. controls (0.5,0.3) and (-0.5,0.3) ..  (-1,0);
\draw[thick, rotate=120] [yshift=0.577cm] (-1,0) .. controls (-1,0.555) and (-0.555,1) .. (0,1)
               .. controls (0.555,1) and (1,0.555) .. (1,0)
               .. controls (0.5,0.3) and (-0.5,0.3) ..  (-1,0);
\draw[thick, rotate=240] [yshift=0.577cm] (-1,0) .. controls (-1,0.555) and (-0.555,1) .. (0,1)
               .. controls (0.555,1) and (1,0.555) .. (1,0)
               .. controls (0.5,0.3) and (-0.5,0.3) ..  (-1,0);
\end{tikzpicture}
}
\caption{}
\label{fig:trefoil}
\end{figure}
This means that the direct analog of the aforementioned Hopf theorem for spherical $f$--neighbors does not hold for small distances.
It would be interesting to find more properties of the set of distances between spherical $f$--neighbors for a continuous map~$f$ of given metric spaces.
For example,  \emph{Is it true that for any continuous map $f\colon \SSS^n \to \R^{n+k}$, the set
$$
\Gamma_f:=\{\delta\in\R\mid \text{$\delta=d(p,q)$ for a pair $\{p,q\}$ of spherical $f$--neighbors}\}
$$
contains a nondegenerate interval?} \emph{Is there a nonzero lower bound for the diameter of $\Gamma_f$?}

\item {\bf Topological Tverberg theorems.} 
Projecting a Euclidean $n$-sphere~$\SSS^n\subset\R^{n+1}$ into a hyperplane in~$\R^{n+1}$ shows that there exists a continuous map $f\colon \SSS^n\to\R^n$ with no set of spherical $f$--neighbors of cardinality exceeding~$2$.
Consequently, each convex $n$-polytope~$P$ has a map $f\colon \dd P\to\R^{n-1}$ with no set of spherical $f$--neighbors intersecting three disjoint faces of~$P$. 
This means that no direct analog of topological Tverberg theorems with three or more disjoint faces holds for spherical $f$--neighbors.
This correlates with the property (mentioned in Remark~\ref{rmrk:principal-disjoint}) that no principal cover has disjoint elements.
Nevertheless, we have some analogs of the topological Radon theorem, which is the topological Tverberg theorem for two disjoint faces: see Corollaries~\ref{cor:cube}--\ref{cor:fullerene}.
It would be interesting to \emph{find extensions of topological Tverberg theorems for spherical $f$--neighbors with additional restrictions}. 
(See also van Kampen--Flores and Conway--Gordon--Sachs type results~\cite{Sko18}.)

\item {\bf Weak Radon rank.} 
{\em Describe the set of polyhedra that are not weakly Radon. Find the weak Radon rank for fullerenes.}

\item {\bf Minimaxes.}
Let $(X,\rho)$ and $(M,d)$ be metric spaces, and let $f\colon X\to M$ be a continuous map. 
Let $P_f$ be the set of all pairs of spherical $f$--neighbors in~$X$.
We set  
$$
D_f:=\sup_{\{x,y\}\in P_f}\rho(x,y),
$$
$$
\mu(X,M):=\inf\limits_{f\in C(X,M)}{D_f}, 
$$
where $C(X,M)$ stands for continuous maps.
Suppose $X=\SSS^n$ and $M=\R^m$. If $m\le n$, then $\mu(\SSS^n,\R^m)=2$ by the Borsuk--Ulam theorem.
For $n<m$, it is shown in~\cite{MM20} that
$$
\mu(\SSS^n,\R^m)= \sqrt{\frac{2(n+2)}{n+1}}. 
$$
It is an interesting problem to {\em find $\mu(X,M)$ and its lower bounds in general and some special cases.} 
In particular, it would be interesting to 
{\em find $D_f$ and~$\mu$ for the case where $M={\mathbb R}^n$ and $X$ is an $n$-dimensional Riemannian manifold.}

\item {\bf Minimaxes 2.} Let us fix $[\mathcal C]$ in $[X,\SSS^{n-2}]$ (see Definition \ref{def2}), for instance, $[\mathcal C]\ne 0$ in $\pi_3(\SSS^2)$. It is an interesting problem, {\em What is min--max distance between 
the points of a set intersecting each element of a cover of this class?} 

\item {\bf Widths, distortion, filling radius, etc.}
Similarly to~$\mu(X,M)$, we consider infima of~$D_f$ over families of homotopic maps, over all continuous maps of a given space to certain classes of spaces (e.\,g., contractible spaces), etc. 
This generates a series of new metric `$\mu$-invariants' of maps and metric spaces.
This $\mu$-invariants are similar to such invariants as distortion, filling radius, various widths, etc. (see \cite{Tikh76,Gro78,Gro83,DS09,Par11,Kar12,AKV12}).
It is an interesting problem {\em to find and describe relations between $\mu$-invariants and classical ones}.



\item {\bf Topological and visual $f$--neighbors.}
Let $f\colon X\to Y$ be a map of topological spaces.
We say that two points $a$ and~$b$ in~$X$ are \emph{topological $f$--neighbors} if $f(a)$ and~$f(b)$ belong to the boundary of the same connected component of the complement~$Y\setminus f(X)$.
If $Y$ is a geodesic metric space, we say that $a$ and~$b$ in~$X$ are \emph{visual $f$--neighbors} if $f(a)$ and~$f(b)$ are connected by a geodesic, in~$Y$, whose interior does not meet~$f(X)$. 
It is interesting to \emph{translate the above constructions and questions to these new types of $f$--neighbors}.

\item {\bf Helly-type sufficient conditions for principal covers.}
Remark~\ref{rmrk:suf-cond} implies some Helly-type sufficient conditions for principal covers.
For example, if $\mathcal{C}=\{C_1,\ldots,C_n\}$ is a closed cover of a normal space~$Y$ such that $\mathcal{N}(\mathcal{C})=\partial\Delta^{n-1}$ and for each $J\subset\{1,\ldots,n\}$ with $|J|\le n-2$ any continuous map $\SSS^{n-2-|J|}\to \bigcap_{j\in J}C_j$ is null--homotopic, then there exists a map $f\colon\partial\Delta^{n-1}\to Y$ such that the image of each facet is contained in an element of $\mathcal{C}$, so that $\mathcal{C}$ is principal. (See the proofs of Theorems~5 and 6 in~\cite{Bog02}.)
It is interesting to find out, \emph{Which of the other versions of topological Helly theorem} (see, e.\,g.,~\cite{Bog02, Mon14}) \emph{give sufficient conditions for principal and non--null--homotopic covers?}
\end{enumerate}

\subsection*{Acknowledgements}
The authors are grateful to Florian Frick, Sergei Ivanov, Roman Karasev, Gaiane Panina, and Arkadiy Skopenkov for helpful discussions and comments.  
Also, the authors are grateful to the anonymous referees for helpful remarks and suggestions.


\begin{thebibliography}{99}

\bibitem{AMS07}
Adams, C., F. Morgan, and J. M. Sullivan,
``When soap bubbles collide.''
\emph{Amer. Math. Monthly} 
{\bf114}:4 (2007), 329--337. 

\bibitem
{ABF20}
Adams, H.,  J. Bush, and F. Frick,  
``Metric thickenings, Borsuk--Ulam theorems, and orbitopes.''
\emph{Mathematika} 
{\bf66}:1 (2020),  79--102. 

\bibitem
{AKV12}
Akopyan, A., R. Karasev, and A. Volovikov,
``Borsuk--Ulam type theorems for metric spaces.''
(2012), preprint arXiv:1209.1249.

\bibitem
{Ale77}
Alexander, R.
``The width and diameter of a simplex.'' 
\textit{Geom. Dedicata} 
{\bf6} (1977), 87--94. 

\bibitem
{BS18}
B\'ar\'any, I. and P. Sober\'on, 
``Tverberg's theorem is 50 years old: a survey.''
\textit{Bull. Amer. Math. Soc. (N.S.)} 
{\bf 55}:4 (2018), 459--492. 

\bibitem
{BBZ16}
B\'ar\'any, I., P. V. M. Blagojevi\'c, and G. M. Ziegler, 
``Tverberg's theorem at 50: extensions and counterexamples.''
\textit{Notices Amer. Math. Soc.} 
{\bf 63}:7 (2016), 732--739.



\bibitem
{Bog02}
Bogatyi, S. A.
``Topological Helly theorem.'' [In Russian.]
\textit{Fundam. Prikl. Mat.}
{\bf8}:2 (2002), 365--405.

\bibitem
{BE18}
Buchstaber, V. M. and N. Yu. Erokhovets, 
``Fullerenes, polytopes and toric topology.'' 
Combinatorial and toric homotopy, 67--178, 
Lect. Notes Ser. Inst. Math. Sci. Natl. Univ. Singap., 35, 
World Sci. Publ., Hackensack, NJ, 2018. 

\bibitem
{DGK63}
Danzer, L., B. Gr\"unbaum, and V. Klee, 
``Helly's theorem and its relatives.'' 
In Proc. Sympos. Pure Math., Vol. VII, pp.~101--180, Amer. Math. Soc., Providence, R.I., 1963.

\bibitem
{Dekster}
Dekster, B. V.
``The Jung theorem for spherical and hyperbolic spaces.''
\textit{Acta Math. Hungar.} {\bf67}:4 (1995), 315--331. 

\bibitem
{DS09}
Denne, E. and J. M. Sullivan,
``The distortion of a knotted curve.''
\textit{Proc. Amer. Math. Soc.}
{\bf137}:3 (2009), 1139--1148. 

\bibitem
{Fri15}
Frick, F.
``Counterexamples to the topological Tverberg conjecture.''
\textit{Oberwolfach Reports} 
{\bf 12} (2015), 318--322.

\bibitem
{Gro78}
Gromov, M. L. 
``Homotopical effects of dilatation.''
\textit{J. Differential Geom.}
{\bf13} (1978), 303--310. 

\bibitem
{Gro83}
Gromov, M. L. 
``Filling Riemannian manifolds.''
\textit{J. Differential Geom.}
{\bf18} (1983), 1--147.


\bibitem
{Gru65}
Gr\"unbaum, B.
``On the facial structure of convex polytopes.''
\textit{Bull. Amer. Math. Soc.} 
{\bf71} (1965), 559--560. 

\bibitem
{Gus}
Gustin, W.
Review of Verblunsky~\cite{Ver52}, MR0051539 (14,495e).

\bibitem
{Iva20}
Ivanov, N. V.
``Leray theorems in bounded cohomology theory.''
(2020), preprint arXiv:2012.08038.

\bibitem
{Kar08}
Karasev, R. N.
``Topological methods in combinatorial geometry.''
\textit{Russian Math. Surveys} 
{\bf 63}:6 (2008), 1031--1078.

\bibitem
{Kar12}
Karasev, R. N.
``A topological central point theorem.''
\textit{Topology Appl.} 
{\bf159}:3 (2012), 864--868.

\bibitem
{Kle}
Klee, V.
Review of Timan~\cite{Tim65}, MR0183000 (32,482).

\bibitem
{Mal16}
Malyutin, A. V.
``On the question of genericity of hyperbolic knots.'' 
\textit{International Mathematics Research Notices}
{\bf2020}:21 (2020), 7792--7828.

\bibitem
{MM20}
Malyutin, A. V. and O. R. Musin, 
``Borsuk--Ulam type theorems for Delaunay maps.'' 
(2020),  manuscript in preparation.

\bibitem
{Mat03}
Matou\v{s}ek, J. 
``Using the Borsuk--Ulam theorem.'' 
Springer-Verlag, Berlin, 2003.

\bibitem
{Mat14}
Matschke, B.
``A survey on the square peg problem.''
\textit{Notices Amer. Math. Soc.} 
{\bf 61}:4 (2014), 346--352. 

\bibitem
{Mon14}
Montejano, L.
``A~new topological Helly theorem and some transversal results.'' 
\textit{Discrete Comput. Geom.} {\bf52}:2 (2014), 390--398. 

\bibitem
{Mor50}
Morita, K.
``On the dimension of normal spaces. II.'' 
\textit{J. Math. Soc. Japan} 
{\bf2} (1950), 16--33.

\bibitem
{Mor75}
Morita, K.
``On generalizations of Borsuk's homotopy extension theorem.''
\textit{Fund. Math.} 
{\bf88}:1 (1975), 1--6. 

\bibitem
{Mun75}
Munkres, J. R.
``Topology: a first course.'' 
Prentice-Hall, Inc., Englewood Cliffs, N.J., 1975. 


\bibitem
{Mus16}
Musin, O. R.
``Homotopy invariants of covers and KKM-type lemmas.''
\textit{Algebr. Geom. Topol.} {\bf16}:3 (2016), 1799--1812. 

\bibitem
{Mus17}
Musin, O. R.
``KKM type theorems with boundary conditions.'' 
\textit{J.~Fixed Point Theory Appl.} {\bf19}:3 (2017), 2037--2049. 

\bibitem
{MW18}
Musin, O. R. and J. Wu,  
``Cobordism classes of maps and covers for spheres.'' 
\textit{Topology and its Applications}  
{\bf 237} (2018), 21--25. 

\bibitem
{Nob20}
Noble, N. 
``Poorly separated infinite normal products''
(2020), preprint arXiv:2002.02483.

\bibitem
{Par11}
Pardon, J.
``On the distortion of knots on embedded surfaces.''
\textit{Ann. of Math. (2)}
{\bf174}:1 (2011), 637--646.

\bibitem
{Sch97}
Schechter, E. 
``Handbook of analysis and its foundations.'' 
Academic Press, Inc., San Diego, CA, 1997.

\bibitem
{Sko18}
Skopenkov, A. B. 
``Topological Tverberg conjecture.'' 
\textit{Russian Math. Surveys} 
{\bf73}:2 (2018), 323--353.

\bibitem
{Sta75}
Starbird, M.
``The Borsuk homotopy extension theorem without the binormality condition.''
\textit{Fund. Math.} 
{\bf87}:3 (1975), 207--211. 

\bibitem
{Ste85}
Steinlein, H. 
``Borsuk's antipodal theorem and its generalizations and applications: A survey.'' 
In A. Granas, editor, M\'{e}thodes topologiques en analyse nonlin\'{e}aire, volume 95 of Colloqu. S\'{e}min. Math. Super., Semin. Sci. OTAN (NATO Advanced Study Institute), pages 166--235, Montr\'{e}al, 1985. Univ. de Montr\'{e}al Press. 

\bibitem
{Ste93}
Steinlein, H. 
``Spheres and symmetry: Borsuk's antipodal theorem.'' 
\textit{Topol. Methods Nonlinear Anal.} 
{\bf1}:1 (1993), 15--33.

\bibitem
{Tikh76}
Tikhomirov, V. M.
``Some Questions of the Approximation Theory.'' [In Russian.]
MSU, Moscow, 1976.

\bibitem
{Tim65}
Timan, T. A.
``Proof of a geometric theorem of Jung, and its analogue in the theory of stochastic processes.'' 
\textit{Uspehi Mat. Nauk} 
{\bf20}:3(123) (1965), 213--218. 

\bibitem
{Ver52}
Verblunsky, S.
``On the circumradius of a bounded set.''
\textit{J. London Math. Soc.} 
{\bf27} (1952), 505--507.

\end{thebibliography}
\end{document}